\def\a             {\alpha}
\def\Ad            {{\mathrm{Ad}}}

\def\be            {\begin{equation}}

\def\bbN           {\mathbb{N}}

\def\bbZ           {\mathbb{Z}}
\def\bfe           {{\bf1}}
\def\can           {\gamma}
\def\canr          {\theta}

\def\cD            {{\mathcal{D}}}
\def\cE            {{\mathcal{E}}}

\def\cZ            {{\mathcal{Z}}}
\newcommand\co[1]  {\overline{{#1}}}

\def\E             {{\mathrm{e}}}
\def\ee            {\end{equation}}
\def\End           {{\mathrm{End}}}
\def\eps           {\varepsilon}

\newcommand\erf[1] {Eq.\ (\ref{#1})}

\def\ext           {{\mathrm{ext}}}
\def\Gtwo          {{\mathrm{G}}_2}
\def\Hom           {{\mathrm{Hom}}}
\def\I             {{\mathrm{i}}}
\def\id            {{\mathrm{id}}}

\def\la            {\lambda}
\def\lan           {\langle}
\def\LG            {{\mathit{LG}}}

\def\LIG           {{\mathit{L}}_I{\mathit{G}}}

\def\Mor           {{\mathrm{Mor}}}

\def\MXN           {{}_M {\cal X}_N}
\def\MXM           {{}_M {\cal X}_M}
\def\MXMa          {{}_M^{} {\cal X}_M^\a}
\def\MXMo          {{}_M^{} {\cal X}_M^0}

\def\MXMp          {{}_M^{} {\cal X}_M^+}
\def\MXMm          {{}_M^{} {\cal X}_M^-}
\def\MXMpm         {{}_M^{} {\cal X}_M^\pm}
\def\MXMmp         {{}_M^{} {\cal X}_M^\mp}

\def\NXN           {{}_N {\cal X}_N}

\def\NXM           {{}_N {\cal X}_M}
\def\NXL           {{}_N {\cal X}_L}
\def\MiXMi         {{}_{M_{1}}{\cal X}_{M_{1}}}
\def\MaXMb         {{}_{M_a}{\cal X}_{M_b}}
\def\MaXMa         {{}_{M_a}{\cal X}_{M_a}}
\def\MbXMb         {{}_{M_b}{\cal X}_{M_b}}
\def\McXMc         {{}_{M_c}{\cal X}_{M_c}}
\def\NXMc         {{}_{N}{\cal X}_{M_c}}
\def\om            {\omega}

\def\ran           {\rangle}

\def\rmA           {{\mathrm{A}}}
\def\rmD           {{\mathrm{D}}}
\def\rmE           {{\mathrm{E}}}
\def\rmv           {{\mathrm{v}}}
\def\rms           {{\mathrm{s}}}
\def\rmc           {{\mathrm{c}}}

\def\sig           {\sigma}

\def\SLZ           {{\mathit{SL}}(2;\bbZ)}
\def\SOf           {{\mathit{SO}}(5)}

\def\SUd           {{\mathit{SU}}(3)}

\def\SUn           {{\mathit{SU}}(n)}
\def\SUN           {{\mathit{SU}}(N)}

\def\SUz           {{\mathit{SU}}(2)}

\def\SUf           {{\mathit{SU}}(4)}

\newcommand\tmat[1]{{}^{{\rm t}} {#1}}
\def\tr            {{\mathrm{tr}}}

\def\tn            {{\tilde{n}}}

\def\Subsection    {Subsection\ }

\def\qed{{\unskip\nobreak\hfil\penalty50
\hskip2em\hbox{}\nobreak\hfil  $\Box$
\parfillskip=0pt \finalhyphendemerits=0\par}\medskip}
\def\proof{\trivlist \item[\hskip \labelsep{\it Proof.\ }]}
\def\endproof{\null\hfill\qed\endtrivlist}

\newcommand\lableq[1]{\label{#1}\end{equation}}
\newcommand\labl[1]{\label{#1}}
\def\typei         {type \nolinebreak I}
\def\typeii        {type \nolinebreak II}
\def\typeiii       {type \nolinebreak III}

\documentclass[11pt]{article}
\usepackage{amssymb,amsfonts,latexsym,epic,eepic}

\begin{document}


\newtheorem{definition}{Definition}[section]
\newtheorem{lemma}[definition]{Lemma}
\newtheorem{corollary}[definition]{Corollary}
\newtheorem{theorem}[definition]{Theorem}
\newtheorem{proposition}[definition]{Proposition}
\newtheorem{conjecture}[definition]{Conjecture}
\newtheorem{assumption}[definition]{Assumption}


\title{Fusion Rules of Modular Invariants\footnote{This contribution is dedicated to
Huzihiro Araki
on the occasion of his seventieth birthday}}

\author{David E.\ Evans\\
\\School of Mathematics\\University of Wales Cardiff\\
PO Box 926, Senghennydd Road\\Cardiff CF24 4YH, Wales, UK}

\date{\today}


\maketitle

|

\begin{abstract}
Modular invariants satisfy remarkable fusion rules. Let
$Z$ be a modular invariant associated to a braided subfactor
 $N\subset M$.  The decomposition of the non-normalized modular invariants
 $Z Z^{*}$ and $Z^{*}Z$ into sums of  normalized modular
invariants is related to the decomposition of  the full induced $M-M$ system
of sectors.
\end{abstract}

\tableofcontents

\section{Introduction}

Suppose $N\subset M$ is a braided type III subfactor; i.e.
say the type III factor $N$ posseses a non-degenerate system $\NXN$
of braided endomorphisms
with the inclusion generated by certain sectors of the system.
Then we know by \cite{R1}  that the system $\NXN$ generates a representation
of the modular group $\SLZ$,
with generators $S$ = \{$S_{\la,\mu} \, ; \,
\la,\mu\in \NXN$ \}, $T$ = \{$T_{\la,\mu} \, ; \,
\la,\mu\in \NXN$ \}. Moreover \cite{BE3, BEK1, E1}
 the inclusion generates a modular
invariant $Z$ through the process of $\alpha$-induction from sectors
of $N$ to sectors of $M$:

\begin{equation}
\label{Zisbb}
Z_{\la,\mu} = \lan \a^+_\la,\a^-_\mu \ran \,, \qquad
\la,\mu\in \NXN \,.
\lableq{Za+a-}

\noindent The right hand side is interpreted as multiplicities
of common sectors in the two inductions, which is clearly thus
a matrix with positive integer entries. It commutes with
both $S$ and $T$ matrices or the representation of
$\SLZ$.  In particular this covers the case of all A-D-E
$\SUz$   modular invariants, and much more besides. We say
that a modular invariant is {\em sufferable} if if can be realised
from a subfactor in this way from $\alpha$-induction
on a braided system of endomorphisms.

Moreover, we can generate from $\alpha$-induction
the sectors $\MXMpm$  which in turn generate the full system
$\MXM$ of $M$. The full system $\MXM$ has
$\sum_{\la, \mu} Z_{\la,\mu}^2$
irreducible sectors and clearly $\alpha$-induction
gives representations of the original $N$-$N$ fusion rules
on $\MXM$. However it is the natural action of the $N$-$N$ sectors
on the corresponding $N$-$M$ sectors $\NXM$ is what gives the
 A-D-E classification and its generalizations. In particular
the trace of Z, $\tr Z$ = $\sum_{\la} Z_{\la,\la}$,
  gives the number of $N$-$M$ sectors in
$\NXM$. Now $\sum_{\la, \mu} Z_{\la,\mu}^2 =
\tr  Z Z^*$.
The matrix $Z Z^*$ is a modular invariant, in that it has
postive integral entries, commutes with representation of
the modular group $\SLZ$, but in general will not be physical
in not  having the vacuum entry normalized to be one.
It is therefore tempting to ask whether we can understand
the full $M$-$M$ system in terms of an analysis of the
modular invariant $Z Z^*$, and an inclusion  $N\subset M_1$
with the full $M$-$M$ system being related to
the chiral $N$-$M_1$ system, just as we understand
the $N$-$M$ system from the modular invariant $Z$.
This was the original motivation in \cite{BER1, BER2}
to write the numerical count $\sum_{\la, \mu} Z_{\la,\mu}^2$
as $\tr{ZZ^*}$, the trace of a modular invariant.

For example for  $\SUz$, where moreover always $Z=Z^*$, we have
the following simple commutative fusion rules for the three
modular invariants at level $16$ labelled by the three Dynkin
diagrams with Coxeter number 18:
\[
Z_{\rmD_{10}}^2=2Z_{\rmD_{10}},
\quad
Z_{\rmD_{10}}Z_{\rmE_7} = Z_{\rmE_{7}}Z_{\rmD_{10}} = 2 Z_{\rmE_{7}^2},
\quad
Z_{\rmE_7}^2=Z_{\rmD_{10}}+Z_{\rmE_7}.
\quad
\]

Consider the subfactor $N\subset M$ describing the
E$_7$ modular invariant.
Indeed the fusion graph  of $\a^+_1$ for the
E$_7$ example has two connected components, D$_{10}$
and E$_7$ with the decomposition
$\tr Z_{\rmE_7}^2= \tr Z_{\rmD_{10}}+ \tr Z_{\rmE_7}$
 reflecting this decomposition of the
$M$-$M$ graph. The aim of this paper is to begin to understand
better the decomposition of the full $M$-$M$ system into its
components via the decomposition of the  matrix $Z Z^*$ into
normalized modular invariants.

If a modular invariant $Z$ is associated to an inclusion
$N\subset M$, then where would we try to understand the
doubled  or modular invariant $Z Z^*$ ? This will  be our
driving principle: For a subfactor $N\subset M$, there
is a natural squaring or iteration procedure $N\subset M \subset M_1$
of the basic construction. Indeed, if the decomposition \erf{Zisbb} formula
for a {\typei}
invariant $Z=B^* B$ in matrix form is related to an inclusion   $N\subset M$,
with dual canonical endomorphism
$\co\iota\iota$,
then it is natural to try to understand the
iteration  $ZZ^*=B^*BB^*B$,
 with the basic construction
$N\subset M \subset M_1$
which has dual canonical endomorphism
$\co\iota\iota\co\iota\iota$.

In the next section we outline our framework of preliminaries in more detail.
In \ Section \ref{iota}  we complete some analysis begun in \cite{BE4} regarding
changing the $\iota$ vertex on the $M$-$N$ graph which will be used for
example in \cite{E1} to understand the Kostant polynomials of \cite{Ko}
from a subfactor point of view. It is not necessary
that a given modular invariant can be realised from a subfactor.
However, even if a modular invariant can be realised from a subfactor it is
not clear what the possible dual canonical endomorphisms are. Nevertheless
there is a simple expression for the sum of all possible dual canonical
endomorphisms in \ Subsection \ref{curious}. This will be used for example in
\cite{EP}
for answering the question of which modular invariants are realisable
in concrete situations with given modular data.
Many subactors can give rise to the same
modular invariant. However, in \ Subsection \ref{general} we consider whether
 sufferable modular invariants can
be realised in canonical ways with natural dual canonical endomorphisms.
We then in \ Section \ref{closer} look at the structure of the products of modular
invariants, in particular $ZZ^*$ and  $Z^*Z$, and how their decomposition into normalised
modular invariants is related to the geometry of the related $M$-$M$
system, its decomposition
into $\MXMpm$ orbits, and the decomposition of $\MXMpm$ into $(ZZ^*)_{0,0}$ and
$(Z^*Z)_{0,0}$ $\MXMo$ orbits respectively.
\ Subsection \ref{sutwo} and \ Subsection \ref{suthree}  contain a discussion
of concrete examples from $\SUz$ and $\SUd$
respectively.
In particular the curious example of the full $M$-$M$ system for the
conformal embedding modular invariant $\SUd_9\subset \mathit(E_6)_1$
where there are six $\MXMp$ orbits in $\MXM$ yet the full system
contains besides three copies of  $\cE^{(12)}_1$, also three copies of the
isospectral graph $\cE^{(12)}_2$.

We can write a sufferable modular invariant $Z$ in terms of rectangular branching matrices
as $Z=B_+^* B_-$ so that $ZZ^*=B_+^*B_-B_-^*B_+$
and $Z^*Z=B_-^*B_+B_+^*B_-$.
We look in \ Subsection \ref{pattern} at the sandwiched $B_\pm B_\pm^*$.
This is a modular invariant for the extended system which is in general not
normalized but its decomposition into normalized modular invariants (usually
permutations) and its relationship to the  decomposition of the full system $\MXM$
into $\MXMpm$ orbits is discussed. Finally in \Subsection \ref{interesting}
we discuss some interesting invariants of $\SUn_n$.
In the conclusions of \cite{BE4} we speculated about
modular invariants which look like \typei\ or \typeii\
but really come from heterotic extensions, i.e.\
for which we have different intermediate local subfactors.
We provide   examples, actually making use
of the heterotic ${\mathit{SO}}(16\ell)_1$ modular
invariants ($\ell=1,2,3,...$) treated in \cite{BE4}, and conformal
inclusions $\SUn_n\subset{\mathit{SO}}(n^2-1)_1$.
The simplest case is
${\mathit{SU}}(7)_7\subset{\mathit{SO}}(48)_1$ and by pulling back the
hetorotic situation on ${\mathit{SO}}(48)_1$ we obtain our strange
heterotic modular invariant on ${\mathit{SU}}(7)_7$ -- which of course must
be symmetric.

\section{Preliminaries}
\labl{prelim}

We cite \cite{EK} as a general reference for operator algebras
and subfactors, and recall the sector setting of \cite{lon2}.
Let $A$ and $B$ be \typeiii\ von Neumann factors.
A unital $\ast$-homomorphism $\rho:A\rightarrow B$
is called a $B$-$A$ morphism. The positive number
$d_\rho=[B:\rho(A)]^{1/2}$ is called the statistical
dimension of $\rho$; here $[B:\rho(A)]$ is the
minimal Jones index \cite{J1} of the subfactor
$\rho(A)\subset B$. If $\rho$ and $\sig$ are $B$-$A$
morphisms with finite statistical dimensions, then
the vector space of intertwiners
\[ \Hom(\rho,\sig)=\{ t\in B: t\rho(a)=\sig(a)t \,,
\,\, a\in A \}  \]
is finite-dimensional, and we denote its dimension by
$\lan\rho,\sig\ran$. Indeed we will only consider morphisms
of finite statistical dimension.
To any $B$-$A$ morphism $\rho$ is assigned a conjugate
 $A$-$B$ morphism $\co\rho$ so that the map  [$\rho$] $\rightarrow$
[$\co\rho$]  is additive, antimultiplicative and idempotent --
generalizing the notion of inversion and conjugate representation
in a group or group dual respectively.

We work with the setting of \cite{BEK1}, i.e.\ we are
working with a \typeiii\ subfactor and
finite system $\NXN\subset\End(N)$ of (possibly degenerately)
braided morphisms which is compatable with the inclusion
$N\subset M$. Then the inclusion is in particular forced to have
finite Jones index and also finite depth (see e.g.\ \cite{EK}).
More precisely, we make the following

\begin{assumption}
{\rm We assume that we have a
\typeiii\ subfactor $N\subset M$
together with a finite system of endomorphisms
$\NXN\subset\End(N)$ in the sense of \cite[Def.\ 2.1]{BEK1}
which is braided in the sense of \cite[Def.\ 2.2]{BEK1}
and such that $\canr=\co\iota\iota\in\Sigma(\NXN)$ for the
injection $M$-$N$ morphism $\iota:N\hookrightarrow M$ and a
conjugate $N$-$M$ morphism $\co\iota$.}
\labl{assbasic}
\end{assumption}

With the braiding $\eps$ on $\NXN$ and its
extension to $\Sigma(\NXN)$ (the set of finite sums of
morphisms in $\NXN$) as in \cite{BEK1}, one can
define the $\a$-induced morphisms $\a^\pm_\la\in\End(M)$
for $\la\in\Sigma(\NXN)$ by the Longo-Rehren formula \cite{LR},
namely by putting
\[ \a_\la^\pm = \co\iota^{\,-1} \circ \Ad
(\eps^\pm(\lambda,\canr)) \circ \lambda \circ \co\iota \,, \]
where $\co\iota$ denotes a conjugate morphism of the
injection map $\iota:N\hookrightarrow M$.
Then $\a^+_\la$ and $\a^-_\la$ extend $\la$, i.e.\
$\a^\pm_\la\circ\iota=\iota\circ\la$, which in turn implies
$d_{\a_\la^\pm}=d_\la$ by the multiplicativity of
the minimal index \cite{L3}.
Moreover, we have $\a_{\la\mu}^\pm=\a_\la^\pm \a_\mu^\pm$
if also $\mu\in\Sigma(\NXN)$, and clearly
$\a_{{\rm{id}}_N}^\pm={{\rm{id}}}_M$.
The morphism  $\a_{\co\la}^\pm$
is a conjugate for $\a_\la^\pm$.
Let $\can=\iota\co\iota$ denote Longo's canonical
endomorphism from $M$ into $N$.

We will assume that braiding on  the system $\NXN$ is non-degenerate. In this
case there is a natural represention of the modular group
$\SLZ$ where the $S$ and $T$ matrices are basically given by the
Hopf link and twist respectively. More precisely,
recall that the statistics phase of $\om_\la$ for
$\la\in\NXN$ is given as
$ d_\la \phi_\la(\eps^+(\la,\la))=\om_\la \bfe $,
where the state $\phi_\la$ is the left inverse of $\lambda$.
We set $z =  \sum_{\la\in\NXN} d_\la^2 \omega_\la$
If $z\neq 0$ we put $c=4\arg(z)/\pi$, which is the
central charge defined
modulo 8.
The $S$-matrix  is defined by
\[ S_{\la,\mu} =  \frac{1}{|z|} \sum_{\rho\in\NXN}
\frac{\om_\la \om_\mu}{\om_\rho} N_{\la,\mu}^\rho d_\rho \,,
\qquad \la,\mu\in\NXN \,,\]
with $N_{\la,\mu}^\rho=\lan\rho,\la\mu\ran$ denoting the
fusion coefficients, \cite{R1,FG1,FRS2}. (As usual, the
label $0$ refers to the identity morphism $\id\in\NXN$.)
Let $T$ be the diagonal matrix with entries
$T_{\la,\mu}=  \E^{-\I\pi c/12 \om_{\la}} \delta_{\la\mu}$.
Then this pair of $S$ and $T$ matrices satisfy
 $TSTST=S$ and give a
unitary representation of the modular group $\SLZ$, 
\cite{R1, T}.
Putting $Z_{\la,\mu}=\lan\a^+_\la,\a^-_\mu\ran$
defines a matrix with positive integral entries
normalized at the vacuum, $Z_{0,0} =1$,
 commuting with $S$ and $T$.
Consequently, $Z$ gives a modular invariant
 \cite{BEK1, E1}.

Let $\MXM\subset\End(M)$ denote a system of endomorphisms
consisting of a choice of representative endomorphisms of
each irreducible subsector of sectors of the form
$[\iota\la\co\iota]$, $\la\in\NXN$.
We choose $\id\in\End(M)$
representing the trivial sector in $\MXM$.
Then we define similarly the chiral systems
$\MXMpm$ and the $\a$-system $\MXMa$ to be the subsystems
of endomorphisms $\beta\in\MXM$
such that $[\beta]$ is a subsector of $[\a^\pm_\la]$ and of
of $[\a_\la^+\a_\mu^-]$, respectively,
for some $\la,\mu\in\NXN$.
The neutral system is defined
as the intersection $\MXMo=\MXMp\cap\MXMm$, so that
$\MXMo \subset \MXMpm \subset \MXMa \subset \MXM$.

Suppose that we have two subfactors, $N \subset M_a$
and $N \subset M_b$ where the irreducible components
of both dual canonical endomorphisms lie in the braided
non degenerate system $\NXN$
with corresponding modular invariants
$Z^a$ and $Z^b$ respectively. Let
$\MaXMb$ denote the irreducible subsectors of
$\iota_a\la\co\iota_b$ where $\iota_a$, $\iota_b$
are the corresponding embeddings of $N$ in $M_a$ and
$M_b$ respectively. We can then by an
extension of the ideas of \cite{BEK1} show that the complexification
of the bimodule $\MaXMb$ under the left action of
$\MaXMa$ and the right action of $\MbXMb$ is isomorphic to
\be
\label{bimoduledec}
\bigoplus_{\la,\mu\in\NXN}
H^a_{\la,\mu} \otimes \co H^b_{\la,\mu},
\ee
where
\be
\label{hilbc}
H^c_{\la,\mu} = \bigoplus_{x\in\NXMc} \Hom(\la\co\mu,x\co x),
\la,\mu\in\NXN.
\ee
is the Hilbert space of intertwiners of dimension $Z_{\la,\mu}^c$, $c = a,b$.
In particular the decomposition in
\erf{bimoduledec} is compatable in the natural way as a bimodule with
 the complexification of the fusion rule algebra of
$\McXMc$ as
\be
\label{fullfudec}
\bigoplus_{\la,\mu\in\NXN} B(H^c_{\la,\mu}).
\ee

\noindent A dimension counts shows that the number of irreducible
$M_a$-$M_b$ sectors of $\MaXMb$ is $\tr$$(Z^{a*}Z^b)$. If $M_a$
=$M_b$ = $M$, then $\#\MXM = \tr Z^*Z$, and if $M_a$ = $N$ and
$M_b$ = $M$, then $\#\NXM = \tr Z$. The action of $\NXN \times
\NXN$ on $\MaXMb$ via $\a$-induction namely $\nu, \rho
\rightarrow \a_\nu^+\a_\rho^-$, on either the left via the
induction $N \subset M_a$ or on the right via $N \subset M_b$,
gives a doubled nimrep $(\nu,\rho) \rightarrow \Gamma_{\nu,
\rho}$ whose spectrum is $S_{\la,\nu}S_{\mu,\rho}/
S_{\la,0}S_{\nu,0}$ with multplicity $Z^a_{\la,\mu}Z^b_{\la,\mu}$.
This reduces to parts 1 and 2 respectively of \cite[Thm.\
4.16]{BEK2} when $M_a = M_b$, $M_a = N$ respectively.
Applications of the existence of such $Z^{a}$-$Z^b$ nimreps for
sufferable invariants and the question of the decomposition of
the products $Z^{a*}Z^b$ into normalised modular invariants  will
appear elsewhere.

We are particularly concerned here with modular invariants
arising in WZW or loop group settings. The modular data ($S$, and $T$
matrices etc) can be
constructed from representation theory of unitary integrable highest
weight modules over affine Lie algebras or in exponentiated form
from the positive energy representations of loop groups.
The subfactor machinery is invoked as follows. Let
 $\LG$ be a loop group (associated to a simple,
simply connected loop group $G$). Let $\LIG$ denote the subgroup
of loops which are trivial off some proper interval $I\subset
S^1$. Then in each level $k$ vacuum representation $\pi_0$ of
$\LG$, we naturally obtain a net of \typeiii\ factors $\{N(I)\}$
indexed by proper intervals $I\subset S^1$ by taking
$N(I)=\pi_0(\LIG)''$ (see \cite{W,FG2,Bg}). Since the
Doplicher-Haag-Roberts DHR selection criterion (cf.\ \cite{H}) is
met in the (level $k$) positive energy representations $\pi_\la$,
there are DHR endomorphisms $\la$ naturally associated with them.
(By some abuse of notation we use the same symbols for labels of
positive energy representations and endomorphisms.) The rational
conformal field theory  RCFT modular data matches that in the
subfactor setting -- in particular
 the RCFT Verlinde fusion coincides with
the (DHR superselection) sector fusion, i.e.\ that
$N_{\la,\mu}^\nu=\langle\la\mu,\nu\rangle$.
The statistics
$S$- and $T$-matrices are identical with the Kac-Peterson
$S$- and $T$-modular matrices which perform the
conformal character transformations.
Antony Wassermann has informed us that
he has extended his results for $\SUn_k$ fusion \cite{W}
to all simple, simply connected loop groups; and with
Toledano-Laredo all but E$_8$ using a variant of the
Dotsenko-Fateev differential equation considered in
his thesis \cite{Tol}, see also \cite{W,Loke,Tol,B1,B2}.

Two subfactor cases are of particular interest in this context,
that of conformal embeddings \cite{X2, BE2, BE3}
and simple current or orbifold constructions \cite{BE2}. For a
conformal embedding $G_k \subset H_1$ we have subfactors $N =
\pi^0({\mathit{L}}_I{\mathit{G}})'' \subset
\pi^0({\mathit{L}}_I{\mathit{H}})'' = M$,  with $\pi^0$ denoting
the level 1 vacuum representation of $\mathit{LH}$. Here, the
subfactor comes equipped with non-degenerately braided systems of
endomorphisms on $N$ and $M$ isomorphic to the level $k$
representations of $G$ and level 1 representations of $H$
respectively, and is relevant for the role of studying conformal
embedding modular invariants. The centre $\bbZ_n$ of $\SUn$ acts
on the algebra $N = \pi_0({\mathit{L}}_I{\SUn})''$, for say the
vacuum level $k$ representation. We can form the crossed product
subfactor $N(I) \subset N(I) \rtimes \bbZ_n$, which will recover
the orbifold modular invariants, but this extended system is only
local
 if and only if $k\in 2n\bbN$
if $n$ is even and $k\in n\bbN$ if $n$ is odd \cite{BE2}.

\section{A closer look at the $M$-$N$ system}
\labl{iota}

For a (non-degenerately) braided subfactor it is the
$M$-$N$ (or $N$-$M$) system which is relevant for the
diagonal part of the modular invariant. Therefore it
is in particular the key to understand the role of
(Coxeter) exponents.

\subsection{Varying the $\iota$-vertex on the $M$-$N$ graphs}

We here assume that we are dealing with a
braided (\typeiii) subfactor $N\subset M$.
For $a\in\NXM$ consider the (irreducible) subfactor
$a(M)\subset N$ and let
\[
a(M) \subset N \subset L
\]
be its basic extension. Note that then $\canr_L=a\co a$
has a Q-system \cite{L4} for $a(M)\subset N$ so that it is a
canonical endomorphism, i.e.\ $\canr_a$ is the
dual canonical endomorphism of $N\subset L$.
Thus $\canr_a=a\co a=\co\iota _L \iota_L$ for
$\iota_L:N\hookrightarrow L$ the injection homomorphism
and $\co\iota _L\in\Mor(L,N)$ a conjugate so that
$\co\iota _L(L)=a(M)$.
We conclude that $\co\iota _L^{-1}\circ a$ is an
isomorphism in $\Mor(M,L)$ with conjugate (i.e.\ inverse)
$a^{-1}\circ\co\iota _L\in\Mor(L,M)$.
For any $b\in\Mor(M,N)$ we now associate $x_b\in\Mor(L,N)$
by putting
\[
x_b=b\circ a^{-1} \circ \co\iota _L.
\]
Note that $x_b$ is irreducible if and only if $b$ is and
that $x_a=\co\iota _L$.

\begin{lemma}
Varying $b\in\NXM$, the $x_b$'s yield all the $N$-$L$ sectors,
and this provides a canonical bijection between $\NXM$ and
$\NXL$.
\end{lemma}

\proof
Note that for any $b\in\NXM$ there is some $\la\in\NXN$
such that $\langle b\co a,\la\rangle\neq 0$ as
$b\co a\in\Sigma(\NXN)$. Thus
\[
\langle x_b,\la\co\iota _L \rangle = \langle b a^{-1}
\co \iota _L \iota_L,\la \rangle = \langle b a^{-1} a \co a,
\la \rangle \neq 0,
\]
implying that $[x_b]$ is one of the $N$-$L$ sectors.
Conversely, assume that there is some $x\in\NXL$ such
that $\langle x,x_b\rangle=0$, i.e.\
$\langle x,ba^{-1}\co\iota _L \rangle=0$ for all $b\in\NXM$.
This implies $\langle x, \la a a^{-1}\co\iota _L\rangle
=\langle x,\la \co\iota _L \rangle \neq 0$ for all $\la\in\NXN$,
in contradiction to $x\in\NXL$.
\endproof

\begin{lemma}
For $b,c\in\NXM$ we have
\be
\langle x_b,\nu x_c \rangle = \langle b,\nu c \rangle,
\ee
i.e.\ the (graphs describing the) multiplication rules of
$\NXN$ on $\NXM$ and $\NXL$ are the same.
\end{lemma}

\proof
This is just
\[
\langle x_b,\nu x_c \rangle = \langle b a^{-1} \co\iota _L,
\nu c a^{-1} \co\iota _L a \rangle =
\langle b, \nu ca^{-1}\co\iota _L \iota_L^{-1}a\rangle
=\langle b,\nu c \rangle,
\]
using that $\iota_L^{-1}a$ is a conjugate morphism of
$a^{-1}\co\iota _L$.
\endproof

Note that the lemma implies in particular that
at least the diagonal part of the coupling matrices
produced from $N\subset M$ and $N\subset L$ are the same.
That in fact the full coupling matrix (and not only
the diagonal part) remains invariant under this
change of the $\iota$-vertex has been shown in \cite{BER2}.

\subsection{A curious identity}
\labl{curious}

We here assume that we are dealing with a non-degenerately
braided ({\typeiii}) subfactor $N\subset M$.
We have seen that for a given braided subfactor $N\subset M$,
realizing a coupling matrix $Z$ and a category of morphisms,
we obtain irreducible subfactors with
dual canonical endomorphisms $\canr_a=a\co a$, $a\in\NXM$,
realizing the same $Z$ \cite{BER2}.
It seems likely to be true that this way we in fact exhaust
all irreducible subfactors producing equivalent categories.

Given a modular invariant matrix $Z$, it is usually not easy
to decide whether it can be realized from a subfactor or not,
and, if yes, how the possible dual canonical endomorphisms
might look like. In the latter case, i.e.\ if there is
some $N\subset M$ realizing $Z$,
at least a statement on the sum of all these endomorphisms
can be made in the following

\begin{proposition}
If the braiding on $\NXN$ is non-degenerate we have the identity
\be
\bigoplus_{a\in\NXM} [a\co a] = \bigoplus_{\la,\mu\in\NXN}
Z_{\la,\mu} [\la\co\mu].
\ee
\end{proposition}

\proof
The multiplicity of $[\nu]$ on the left-hand side is
for all $\nu\in\NXN$
\[
\sum_{a} \langle a\co a,\nu \rangle
= \sum_{a} \langle a,\nu a \rangle
= \tr (G_\nu) = \sum_\rho Z_{\rho,\rho}
\frac{S_{\rho,\nu}}{S_{\rho,0}},
\]
where we used \cite[Thm.\ 4.16]{BEK2}.
The multiplicity of $[\nu]$ on the right-hand side is
for all $\nu\in\NXN$
\[
\sum_{\la,\mu} Z_{\la,\mu}
\langle \la\co\mu,\nu \rangle =
\sum_{\la,\mu} Z_{\la,\mu} N_{\nu,\mu}^\la
= \sum_{\la,\mu,\rho} Z_{\la,\mu}
\frac{S_{\rho,\nu}}{S_{\rho,0}} S_{\rho,\mu} S_{\rho,\la}^*
= \sum_{\rho} Z_{\rho,\rho}
\frac{S_{\rho,\nu}}{S_{\rho,0}},
\]
where we used the Verlinde formula and modular invariance.
\endproof

\subsection{Towards a general formula for $[\canr]$}
\labl{general}

Looking at a couple of examples, it seems that a ``physical
invariant'' $Z$, which can be realized from some subfactor,
can in fact be realized with a dual canonical endomorphism
given by something like
\be
[\canr] = \bigoplus_\la Z_{\la,\co\la} [\la].
\lableq{strange}
In general summing over a subset
of $\NXN$ related to Frobenius-Schur indicators and
conformal dimensions.  Let us consider some examples.

For $\bbZ_n$ conformal field theories with $n$ odd, this works perfectly.
In this situation, there are
$n$ sectors, labelled by $\la=0,1,2,...,n-1$ (mod$\,n$),
obeying $\bbZ_n$ fusion rules, and conformal dimensions
of the form $h_\la=a\la^2/2n$ (mod$\,1$), where $a$ is an integer
mod$\,2n$, $a$ and $n$ coprime and $a$ is even whenever
$n$ is odd. The modular invariants of such models have been
classified \cite{D}. They are labelled
(with notation as in \cite{BER1,BER2}) by the divisors
$\delta$ of $\tn$, where $\tn=n$ if $n$ is odd and $\tn=n/2$ if
$n$ is even. Let us take $n$ odd.
Then it is not hard to show that
$Z^{(\delta)}_{\la,\co\la}=Z^{(n/\delta)}_{\la,\la}$.
Thus by \cite[Eq.\ (8.1)]{BER1}  we find
$Z^{(\delta)}_{\la,\co\la}=1$ for $\la=0$ mod $\delta$
and $Z^{(\delta)}_{\la,\co\la}=0$ otherwise, and in fact
$\canr=\bigoplus_{j=0}^{n/\delta-1}[\rho_{j\delta}]$
realizes $Z^{(\delta)}$, see \cite{BER1}.
(By the way: Since
$\langle\a^+_j\a^-_{-j},\can\rangle=
\langle\a^+_j\a^+_{-j},\can\rangle=1$
it is easy to see that
$[\can]=\sum_{j=0}^{\tilde n /\delta -1} [\a^+_j\a^-_{-j}]$
for all $\bbZ_n$ theories, no matter whether $n$ is
even or odd.)

Note that for the conjugation invariant $C$ we would
usually insert  all morphisms in the $[\canr]$. This does not
work for the $\bbZ_n$ CFT's with $n$ even because we must
use the even labels only \cite{BER1}. So for some
reasons the odd labels have to be ruled out.
(Moreover, if $n$ is a multiple of 4 we do not  want to
see the self-conjugate even label $n/2$ in the dual
canonical endomorphism realizing the trivial invariant.)
A similar thing happens for $\SUz_k$.
Here $Z_{\la,\la}=Z_{\la,\co\la}$,
but if we restrict the sum to even spins then we can
in fact realize each A-D-E invariant by \erf{strange}.

Let us start with the subfactors used to produce the A-D-E modular
invariants in \cite{BE2,BEK1,BEK2}, i.e.\ the corresponding
dual canonical endomorphisms $[\canr]$ are given by:
\[
\begin{array}{lll}
\mathrm{A}_{\ell}, \qquad & \ell=k+1 \, : \qquad & [\la_0] \\
\mathrm{D}_{\ell}, & k= 2 \ell -4 \, : & [\la_0] \oplus [\la_k] \\
\mathrm{E}_6, & k= 10 \, : & [\la_0] \oplus [\la_6] \\
\mathrm{E}_7, & k= 16 \, : & [\la_0] \oplus [\la_8] \oplus [\la_{16}] \\
\mathrm{E}_8, & k= 28 \, : & [\la_0] \oplus [\la_{10}] \oplus
                      [\la_{18}] \oplus [\la_{28}]
\end{array}
\]
Now we choose the following $M$-$N$ morphisms $[\co a]$:
For A$_\ell$ (where $\iota$ is trivial) we choose
$\iota\la_{[k/2]}\equiv\la_{[k/2]}$.
Here $[x]$ denotes the greatest possible
integer less than or equal to $x$.
For D$_\ell$ we choose $\iota\la_{[\ell/2]-1}$.
For E$_6$ we choose $\sigma\iota$ with $\sigma$
the marked vertex with statistical dimension $\sqrt2$.
For E$_7$ we choose the morphism denoted by $\co b '$
in \cite[Fig.\ 41]{BEK2}.
For E$_8$ we choose $\a^{(1)}_6\iota$ with
$\a^{(1)}_6$ the neutral or marked vertex as in \cite[Fig.\ 8]{BE2}.
It is now straightforward to compute the sectors
$[a\co a]$ which will be our new $[\canr]$'s.
For example, for D$_\ell$ we compute
$[\la_{[\ell/2]-1}\co\iota\iota\la_{[\ell/2]-1}]
=[\la_{[\ell/2]-1}]^2([\la_0]\oplus[\la_k])$.
For E$_6$ we compute
$[\co\iota\sigma\sigma\iota]=[\co\iota]([\a_0]\oplus[\a_{10}])
[\iota]=([\la_0]\oplus[\la_{10}])([\la_0]\oplus[\la_6])$.
Only for E$_7$ we need to sit down a bit, using
$[\co b ']=[\iota\la_2]\oplus[\iota\la_4]\ominus[\iota\la_6]$.
This gives:
\[
\begin{array}{lll}
\mathrm{A}_{\ell}, & \ell=k+1 \, : &
 [\la_0] \oplus
 [\la_2] \oplus [\la_4] \oplus \ldots \oplus [\la_{2[k/2]}] \\
\mathrm{D}_{2\varrho}, & k= 4 \varrho -4 \, : & [\la_0] \oplus
 [\la_2] \oplus \ldots \oplus [\la_{2\varrho-4}]\oplus
 2[\la_{2\varrho-2}] \oplus [\la_{2\varrho}] \oplus
 \ldots \oplus [\la_k] \\
\mathrm{D}_{2\varrho+1}, & k= 4 \varrho -2 \, : & [\la_0] \oplus
 [\la_2] \oplus [\la_4] \oplus \ldots \oplus [\la_k] \\
\mathrm{E}_6, & k= 10 \, : & [\la_0] \oplus [\la_{4}] \oplus
 [\la_6] \oplus [\la_{10}] \\
\mathrm{E}_7, & k= 16 \, : & [\la_0] \oplus [\la_4] \oplus
 [\la_6] \oplus [\la_8] \oplus [\la_{10}] \oplus
 [\la_{12}] \oplus [\la_{16}] \\
\mathrm{E}_8, & k= 28 \, : & [\la_0] \oplus [\la_6] \oplus
 [\la_{10}] \oplus [\la_{12}] \oplus [\la_{16}] \oplus
 [\la_{18}] \oplus [\la_{22}] \oplus
[\la_{28}] .
\end{array}
\]
So here we indeed find exactly the even spins of the diagonal.
(Note that the $[\canr]$'s for A and D$_\mathrm{odd}$ are
the same (at levels $k=6,10,14,...$). Thus these are examples
for subfactors producing different $Z$'s but having the
same dual canonical endomorphism sector.)

It is interesting to note what the canonical endomorphism looks like
in these possibly natural subfactors:

\[
\begin{array}{lll}
\mathrm{A}_{\ell}, & \ell=k+1 \, : &
 [\a_0] \oplus
 [\a_2] \oplus [\a_4] \oplus \ldots \oplus [\a_{2[k/2]}] \\
\mathrm{D}_{2\varrho}, &k=4\varrho -4 \, : & [\a_0] \oplus
 [\a_2] \oplus \ldots \oplus [\a_{2\varrho-4}]\oplus
 [\a_{2\varrho-2}^{(1)}] \oplus [\a_{2\varrho-2}^{(2)}]
 \oplus [\epsilon] \oplus \\
&& \oplus[\beta_2] \oplus [\beta_4] \oplus \ldots \oplus
[\beta_{2\varrho-4}]\oplus
 [\eta] \oplus [\eta'] \\
\mathrm{D}_{2\varrho+1}, & k= 4 \varrho -2 \, : & [\a_0] \oplus
 [\a_2] \oplus [\a_4] \oplus \ldots \oplus [\a_k] \\
\mathrm{E}_6, & k= 10 \, : & [\a_0] \oplus [\a_{10}] \oplus
 [\delta] \oplus [\delta'] \\
\mathrm{E}_7, & k= 16 \, : & [\a_0] \oplus [\eta] \oplus
 [\delta] \oplus [\a_3^+\a_1^-] \oplus [\a_{6}^+] \oplus
 [\a_{4}^+] \oplus \\
 &&\oplus([\a_5^+\a_1^-] \ominus [\a_3^+\a_1^-])\\
\mathrm{E}_8, & k= 28 \, : & [\a_0] \oplus [\a_6^{(1)}] \oplus
 [\delta] \oplus [\chi] \oplus [\omega] \oplus
 [\varpi] \oplus [\eta] \oplus
[\eta'] .
\end{array}
\]

\noindent For $\mathrm{D}_{2\varrho}, \mathrm{D}_{2\varrho+1},
\mathrm{E}_6,\mathrm{E}_7,\mathrm{E}_8$,
we have used the notation of
\cite[Fig. 9]{BE3}, \cite[Fig.40]{BEK2}, \cite[Fig.2]{BE3}, \cite[Fig. 42]{BEK2}
\cite[Fig. 5]{BE3} respectively.

At least in this $\SUz$ setting,
there is a fusion rule
symmetry on $\MXM$
obtained by interchanging $\a^+_\la$  with $\a^-_\la$ taking
$\MXMp$ to $\MXMm$. In terms of the above figures
for $\mathrm{D}_{2\varrho}, \mathrm{D}_{2\varrho+1},
\mathrm{E}_6,\mathrm{E}_7,\mathrm{E}_8$, this is the flip around the
vertical through the vacuum. (When we change in the above examples the
subfactor $N \subset M$ but retain the same modular invariant the systems of
sectors $\MXMpm, \MXMo, \MXM, \MXN$ remain isomporhic to the old ones, so
we retain the same figures). Again for $\mathrm{E}_7$ we need to do some
work, e.g. $\langle \co b' b', \a^+_i\a^-_j \rangle =
 \langle \iota([\la_2] \oplus [\la_4] \ominus [\la_6])^2 \co\iota, \a^+_i\a^-_j \rangle
=  \langle ([\la_2] \oplus [\la_4] \ominus [\la_6])^2, \co\iota
\a^+_i\a^-_j \iota \rangle = \langle ([\la_2] \oplus [\la_4] \ominus
[\la_6])^2, [\la_i][\la_j]([\la_0] \oplus [\la_8] \oplus [\la_{16}]) \rangle$.
Then the
the "real" part of the full system are the sectors
fixed under the flip, i.e. the sectors lying on the vertical through
the vacuum.
Then the canonical endomorphism
is the even part of the "real" part of the full system
- presumably these are the ones of Frobenius-Schur indicator one.

For $\SUd_k$ we have checked for $k=1$ and $k=2$ that
$Z=C$ is indeed realized by \erf{strange}, the sum taken
over all $\SUd_k$ weights. Assuming that the subfactor
exists for $k=3$ it is easy to check it for this case as
well. Similarly is is easy to check for $k\le3$ that
that $[\canr]$ given as sum over all selfconjugate sectors
indeed produces $Z=\bfe$ since the $[\canr]$ is just the
square of the only non-trivial self-conjugate sector
$[\la_{(2,1)}]$. Squaring larger $[\la]$'s instead,
this procedure should also work at any higher levels $k$.

\section{A closer look at the $M$-$M$ system}
\labl{closer}

Here we discuss the structure of the entire $M$-$M$ system.
We will only consider proper modular invariants here,
i.e.\  we assume the $N$-$N$ system is non-degenerate.
First some observations.

\subsection{Some remarks on products of modular invariants}

A modular invariant from a subfactor is of the form
\[
Z_{\la,\mu} = \sum_{\tau\in\MXMo} b^+_{\tau,\la}
b^-_{\tau,\mu}.
\]
Now let consider the fusion graph of $\a^+_\la$ in the
entire system $\MXM$. (We will here consider the non-degenerate
case only.) We know since \cite{BEK2} that the multiplicity
of the eigenvalue $S_{\la,\rho}/S_{0,\rho}$ is given by
$\sum_\mu Z_{\rho,\mu}^2=(ZZ^*)_{\rho,\rho}$, and that this
exhausts the spectrum. Since this graph contains the chiral
graph as a subgraph, we must have
$(ZZ^*)_{\rho,\rho}\ge Z^+_{\rho,\rho}$, where $Z^+$ denotes
the \typei\ parent,
$Z^+_{\la,\mu}=\sum_\tau b^+_{\tau,\la}b^+_{\tau,\mu}$.
And indeed, we can compute quite generally
\[
\begin{array}{rl}
(ZZ^*)_{\la,\mu} &= \sum_\nu Z_{\la,\nu} Z_{\mu,\nu}
= \sum_\nu \sum_{\tau,\tau'} b^+_{\tau,\la}
b^-_{\tau,\nu} b^+_{\tau',\mu} b^-_{\tau',\nu} \\[.4em]
& \ge \sum_\nu \sum_\tau b^+_{\tau,\la}
b^-_{\tau,\nu} b^+_{\tau,\mu} b^-_{\tau,\nu}
 \ge \sum_\tau b^+_{\tau,\la} b^+_{\tau,\mu}
= Z^+_{\la,\mu},
\end{array}
\]
where we used that for each $\tau\in\MXMo$ there is a $\nu$
such that $b^-_{\tau,\nu}\ge1$.
(And of course we obtain similarly
$(Z^*Z)_{\la,\mu}\ge Z^-_{\la,\mu}$, etc.\ etc.)
Note that $ZZ^*-Z^+$ must be modular invariant, and looking
at the above calculation we see that it is even non-negative.
So what about normalization? We distinguish the two cases:
(1) $Z$ is a pure permutation. Then in fact $ZZ^*=Z^+=\bfe$.
(2) otherwise there is a $\la$ with $Z_{0,\la}\neq0$,
and consequently $(ZZ^*)_{0,0}=\sum_\la Z_{0,\la}^2>1$.
If this gives exactly 2 then we know that $ZZ^*-Z^+$
is another normalized integral modular invariant, but if it is
larger then it is not clear  whether $ZZ^*-Z^+$ can always
written as a positive integer linear combination of
normalized integral modular invariants.

It is clear that if we always obtain a
positive integer linear combination of normalized
integral modular invariants, then the number of
such invariants (counting multiplicities) will be
$(ZZ^*)_{0,0}=\sum_\la Z_{0,\la}^2=\langle\canr_+,\canr_+\rangle$.
Each invariant is expected to correspond to a component
of the full fusion graph, so we expect
$\sum_\la Z_{0,\la}^2$ components.
(By components we  mean here a connected component of
the fusion graph of a generator $\a^+_f$. Equivalently
one can decompose the sum of the full fusion matrices
of all chiral sectors into irreducible components.)
That at least this numbering for the connected components
is indeed correct is shown in the following subsection.

\subsection{On the geometry of the $M$-$M$ system}

Let us recall that the $M$-$M$ system has subsystems
$\MXM\supset\MXMpm\supset\MXMo$. Under the action
(fusion) of a chiral system, say $\MXMp$, the $\MXM$
system decomposes into $\MXMp$ orbits. These correspond
to the connected components of the fusion graph
of a generator of $\MXMp$ in $\MXM$. We may draw
such a graph using straight lines, and the graph
arising from the corresponding generator of $\MXMm$
using dotted lines as in \cite[Figs.\ 2,5,8,9]{BE3}
or \cite[Figs.\ 40,42,43]{BEK2}.
For the E$_8$ example we find  4 $\MXMp$
orbits which are precisely the 4 straight-lined
E$_8$ ``layers'' in \cite[Fig.\ 5]{BE3}.
How many such layers do we usually have?
A first answer is this:

\begin{lemma}
The number of $\MXMpm$ orbits in $\MXM$ is equal
to the number of $\MXMo$ orbits in $\MXMmp$.
In fact, all $\MXMpm$ orbits in $\MXM$ intersect
with the subset $\MXMmp\subset\MXM$, and the intersections
are precisely the $\MXMo$ orbits in $\MXMmp$.
\end{lemma}

\proof
Consider the identity component $\Gamma^+_{(0)}$ of the fusion
graph of $\MXMp$ in $\MXM$ (which is essentially $\MXMp$ itself).
Since $\MXMp$ and $\MXMm$ generate $\MXM$, each
connected component $\Gamma^-_{(j)}$ of the fusion graph
of $\MXMm$ in $\MXM$ must touch $\Gamma^+_{(0)}$ somewhere.
(E.g.\ the identity component $\Gamma^-_{(0)}$ meets
$\Gamma^+_{(0)}$ exactly on the ambichiral vertices.)
Hence the number of $\MXMm$ orbits in $\MXM$ is
equal to the number of groups groups of vertices on
$\Gamma^+_{(0)}$ lying on the same component $\Gamma^-_{(j)}$.
Two vertices on $\Gamma^+_{(0)}$ corresponding to
sectors $\beta_1,\beta_2\in\MXMp$ lie on the same
component $\Gamma^-_{(j)}$ if and only if there is
a $\beta\in\MXMm$ such that
$\langle\beta_1\beta,\beta_2\rangle\neq0$.
But $\langle\beta,\co{\beta_1}\beta_2\rangle\neq0$
if means that $\beta$ is ambichiral.
Hence two vertices on $\Gamma^+_{(0)}$ corresponding to
sectors $\beta_1,\beta_2\in\MXMp$ lie on the same
component $\Gamma^-_{(j)}$ if and only if they are in
the same ambichiral orbit. The proof is completed by
exchanging $+$ and $-$ signs.
\endproof

A more concrete answer is now obtained in the following

\begin{lemma}
The number of $\MXMo$ orbits in $\MXMpm$
is given by $\sum_\la (b^\pm_{0,\la})^2$.
\end{lemma}

\proof
Let $\Gamma^{\pm}_{\tau,0}$ be the fusion matrix
of $\tau\in\MXMo$ in $\MXMpm$, as in \cite[Sect.\ 4]{BEK2}.
The sum matrix $Q=\sum_\tau \Gamma^{\pm}_{\tau,0}$ will
not be irreducible as long as we have more than one
$\MXMo$ fusion orbit (i.e.\ as long as $\MXMo\neq\MXMpm$).
In fact $Q$ must decompose into a number of irreducible
blocks which is exactly the number of fusion orbits.
Nevertheless the vector $\vec{d}$ with entries $d_\beta$,
$\beta\in\MXMp$ is an eigenvector of $Q$ with eigenvalue
$\sum_\tau d_\tau$. Since all the entries are strictly positive,
it must be the direct sum of the Perron-Frobenius eigenvectors
of each irreducible block (up to a scaling by a positive factor
for each block). Thanks to the Perron-Frobenius theorem,
the number $\sum_\tau d_\tau$ is thus the (non-degenerate)
Perron-Frobenius eigenvalue of each irreducible block.
It follows that the number of irreducible components is
given by the multiplicity of the eigenvalue $\sum_\tau d_\tau$,
i.e.\ by the multiplicity of $\chi^\ext_0(\tau)$ in
$\Gamma^{\pm}_{\tau,0}$, $\tau\in\MXMo$.
By the diagonalization of the $\Gamma^{\pm}_{\tau,0}$'s
derived in \cite[Thm.\ 4.16]{BEK2},
we know that this multiplicity  is exactly
$\sum_\la (b^\pm_{0,\la})^2$.
\endproof

Note that $\sum_\la (b^+_{0,\la})^2=\sum_\la Z_{\la,0}^2$
and $\sum_\la (b^-_{0,\la})^2=\sum_\la Z_{0,\la}^2$.
In fact, the consideration of the $\MXMo$ in $\MXMpm$
was instructive, but not really necessary to get the
number of $\MXMpm$ orbits in $\MXM$.
Thanks to the generating property, we could
also have determined the number of $\NXN$ orbits
in $\MXM$ via the induced $[\a^+_\la]$and $[\a^-_\la]$.
Then the statement of \cite[Thm.\ 4.14]{BEK2}
would similarly determine the multiplicity of
the Perron-Frobenius eigenvalues $\chi_0(\la)$
as $\sum_\mu Z_{0,\mu}^2$ and $\sum_\mu Z_{\mu,0}^2$,
respectively.

\section{Examples}

Suppose  $N\subset M$ is a braided subfactor with
$\iota:N\hookrightarrow M$ being the injection map and basic construction
  $N\subset M \subset M_1$. Thus if $\iota_1:M\hookrightarrow M_1$ is the
corresponding inclusion,
then by naturality  the sector $[\co\iota_1\iota_1]$ of
the dual canonical endomorphism for  $M\subset M_1$
is identified with the sector of the canonical endomorphism for
$N \subset M$,
 i.e. $\iota\co\iota$. Hence
the sector of the dual canonical endomorphism $\theta_1$ for
$N \subset M_1$  is
$[\co\iota{\co\iota_1}\iota_1\iota]$
= $[\co\iota\iota\co\iota\iota]$ = $[\theta^2]$,
which lies in $\Sigma(\NXN)$ as $\theta$ does. In particular,
if $\NXN$ is braided, we can certainly apply $\alpha$-induction
to the inclusion
$N \subset M_1$.
Note that in this context, that the inclusion, $N  \subset M_1$
rarely satisfies chiral locality by Corollary 3.6 \cite{BE1}.
We have the
naturality equations for $\a$-induced morphisms
\[
x \eps^\pm(\rho,\la)=\eps^\pm(\rho,\mu)\a^\pm_\rho(x)
\]
whenever $x\in\Hom(\iota\la,\iota\mu)$
and $\rho\in\Sigma(\NXN)$, see e.g
\cite[Eq.\ (9)]{BE4}.
In particular, inducing from $N$ to $M_1$, we have taking
$\la = \mu = id$, and  $x\in\Hom(\iota_1\iota,\iota_1\iota)$, that
$\alpha_{{\rho}}(x) = x$ on $N'\cap M_1$, for all $\rho$.

We will look again at the ${\SUd}$ and ${\SUz}$ situations in detail
in this basic construction.

\subsection{${\SUz}$-invariants}
\labl{sutwo}

By the A-D-E classification \cite{CIZ2}, we know that there are at
most three invariants for each level labelled by
Dynkin diagrams. They satisfy the following fusion rules:

\[
Z_{\rmD_{2\varrho}}^2=2Z_{\rmD_{2\varrho}},
\quad
Z_{\rmD_{2\varrho+1}}^2=Z_{\rmA_{4\varrho-1}},
\quad
Z_{\rmE_6}^2=2Z_{\rmE_6},
\quad
Z_{\rmE_7}^2=Z_{\rmD_{10}}+Z_{\rmE_7},
\quad
Z_{\rmE_8}^2=4Z_{\rmE_8}.
\]

\noindent (i) {\it Example:} D$_j$

 We start with  ${\SUz}$ at even  level $k$  and   the simple current or
orbifold invariants.  Here there is
a  $\bbZ_2$ extension: $N \subset N \rtimes\bbZ_2$,
with $N$ as $ \pi^0({\mathit{L}}_I{\mathit{SU}}(2))''$
in the vacuum representation at level $k$,
and dual canonical endomorphism
$[\la_0] \oplus [\la_k]$.
If $k = 4l -4$, \cite{BE2}  then the extension is local, the corresponding
modular invariant is  D$_{2\ell}$, and the canonical endomorphism
is
$\gamma = [id] \oplus [\alpha^{\pm}_k]$.
If $k = 4l -2$, \cite{BE2}  then the extension is not local, the corresponding
modular invariant is  D$_{2\ell + 1}$ and the
canonical endomorphism  is
$\gamma = [id] \oplus [\epsilon]$.
where $\epsilon$ is an irreducible subsector of
${[\alpha^+_{{1}} \alpha^-_{{1}}]}$.
In either case, the basic construction is by Takesaki duality:
\[
N \subset N \rtimes\bbZ_2 \subset  N \rtimes\bbZ_{2} \rtimes {\hat\bbZ}_{2}  =
{N \otimes  Mat_2}.
\]

Thus by the above naturality, $\alpha_\la$ = $\la \otimes id$,
as here  $ N'\cap M_1$ =  Mat$_2$, the 2 $\times$  2 complex matrices.
Thus $\NXN$  is identified with $\MiXMi$, and we do not appear
to have anything interesting. To see the finer structure,
we need to look closer at the dual canonical endomorphism
$[\theta_1]$, which decomposes in the local case
$k = 4l -4$, into
$[\co\iota\iota]$ and $[\co\iota\alpha^{\pm}_{k}\iota]$. Both are
dual canonical endomorphisms in their own right. The first
can be thought of as giving the sheet of D$_j$ in the full
$\MXM$ system starting at $[\id_M]$
and the second sector as giving the other sheet in
the full $\MXM$ system. All this becomes clearer
in the \typei\ conformal embedding modular invariants.

\noindent (ii) {\it Example:} E$_6$,  $\SUz_{10}\subset\SOf_1$

We now consider the  E$_6$  modular invariant for
 $\SUz$:
\[
Z_{\mathrm{E}_6}= |\chi_0 + \chi_6|^2 + |\chi_4 + \chi_{10}|^2
+ |\chi_3 + \chi_7|^2 \,.
\]

 This is exhibited by the
 conformal embedding
$\SUz_{10}\subset\SOf_1$.
Here the dual canonical endomorphism $\theta$ is given
by the vacuum sector
$[\theta] = [\la_0] \oplus [\la_6]$,
and the corresponding
canonical endomorphism was computed in \cite{BE3}
as
$[\can]=[\id]\oplus[\alpha^+_{{1}}\alpha^-_{{1}}]$.
 Then for the corresponding basic construction
 $N\subset M \subset M_1$
we have

\[
[\theta_1]=
[\co\iota{\co\iota_1}~ \iota_1\iota] = [\co\iota\iota]\oplus\ {[\co\iota\alpha^-_{{1}}\alpha^-_{{1}}\iota]}.
\]

This time the dual canonical endomorphism $[\co\iota\iota]$ gives the first sheet of
the full $\MXM$ system, whilst the second term gives the second sheet of the full system
where the sector $[\alpha^+_{{1}}\alpha^-_{{1}}]$ in the full system is identified
with the $N$-$M$ sector $[\alpha^-_{{1}}\iota]$ using the changing the $\iota$
vertex argument.

\noindent (iii) {\it Example:} E$_8$,  $\SUz_{28}\subset (\Gtwo)_1$

Next let us revisit the E$_8$  modular
invariant at level $k$ = 28:

\[ Z_{\mathrm{E}_8}= |\chi_0 + \chi_{10} + \chi_{18} + \chi_{28}|^2
+ |\chi_6 + \chi_{12} + \chi_{16} + \chi_{22}|^2 \,. \]

\noindent This is exhibited by the conformal embedding
$\SUz_{28}\subset (\Gtwo)_1$.
The dual canonical endomorphism is again given by the vacuum sector

\[ [\theta] = [ \la_0] \oplus[ \la_{10}] \oplus[ \la_{18}] \oplus [\la_ {28}]\,. \]

\noindent The canonical endomorphism was computed in \cite{BE3}
as:

\[
[\can]=[\id_M]\oplus [\alpha^+_{{1}}\alpha^-_{{1}}]
\oplus [\alpha^+_{{2}}\alpha^-_{{2}}]
 \oplus[\eta].
\]

\noindent where [$\eta$] was described as  an irreducible subsector of
$[\alpha^+_{{3}}\alpha^-_{{3}}]$. However by comparing
Fig. 8 of \cite{BE2}
  with Fig. 5 of \cite{BE3}, and with the above experience for
$E_6$, we suspect that $[\eta]$ can be identified
with
$[\alpha_{5}^{+(2)}\alpha_{5}^{-(2)}]$,
where
$[{\alpha_{{5}}^{\pm}}] =[{\alpha_{{5}}^{{\pm}({1})}}] \oplus
[\alpha_{5}^{\pm(2)}]$ and
$[{\alpha_{{7}}^{\pm}}] =
[{\alpha_{{3}}^{{\pm}}}] \oplus
[{\alpha_{{5}}^{{\pm}({1})}}]$.
 We can compute
\[
\begin{array}{rl}
\langle \gamma,\alpha_{5}^{+(2)}\alpha_{5}^{-(2)}\rangle
&= \langle  \iota\co\iota, \alpha_{5}^{+(2)} \alpha_{5}^{-(2)} \rangle\\
&= \langle  id, \co\iota\alpha_{5}^{+(2)}\alpha_{5}^{-(2)}\iota \rangle \\
&=\langle id, \co\iota [\alpha^+_{5} \oplus \alpha^+_{3} \ominus \alpha^+_{{7}}] [\alpha^-_{{5}} \oplus \alpha^-_{{3}} \ominus \alpha^-_{{7}} ] \iota \rangle \\
&= \langle  id, \co\iota\iota([\la_5] \oplus [\la_3] \ominus [\la_7])^2
\rangle\\
&= 1
\end{array}
\]

\noindent using the Verlinde fusion rules for ${\SUz}$ at level 28.
We can similarly show that $[\alpha_{5}^{+(2)}\alpha_{5}^{-(2)}]$
is irreducible and disjoint from $ [\id_M], [\alpha^+_{{1}}\alpha^-_{{1}}]$ and
$ [\alpha_2^+\alpha_2^-]$. Hence 
$[\eta]$=$[\alpha_{5}^{+(2)}\alpha_{5}^{-(2)}]$.

There are four sheets in the full system $\MXM$, all copies of E$_8$.
The four terms in $\gamma$ give rise to the four sheets in the full
system with vertices
$[\id_M]$, $[\alpha^+_{1}\alpha^-_{1}]$,
$[\alpha^+_{{2}}\alpha^-_{{2}}]$,
$ [\alpha_{5}^{+(2)}\alpha_{5}^{-(2)}]$
identified with base points
$\iota$,      $ \alpha^-_{{1}}\iota$,
  $ \alpha^-_{{2}}\iota$ and  ${\alpha_{{5}}^{-({2})}}\iota$
 on the $N$-$M$ graph E$_8$ using again the
argument of changing the $\iota$ vertex.

\subsection{${\SUd}$-invariants}
\labl{suthree}

We now move on the the case of ${\SUd}$ and its modular invariants.

\noindent (i) {\it Example:} ${\cE^{(8)}}$,  $\SUd_5\subset \mathit{SU}(6)_1$.

The first conformal embedding invariant is at level 8:
$$
\begin{array}{rl}
Z_{\mathcal{E}^{(8)}} &= |\chi_{0,0} + \chi_{4,2}|^2 +
|\chi_{ 2,0} + \chi_{5,3}|^2 + |\chi_{2,2} + \chi_{ 5,2}|^2
 + \,|\chi_{3,0} + \chi_{ 3,3}|^2 +\\
& |\chi_{ 3,1} + \chi_{5,5}|^2 + |\chi_{ 3,2} + \chi_ {5,0}|^2
\end{array}
$$

This can be obtained from the conformal inclusion
$\SUd_5 \subset \mathit{SU}(6)_1$
with  dual canonical endomorphism is given by the vacuum sector
$[\theta] = [\la_{0,0}] \oplus [\la_{4,2}]$
 with the canonical endomorphism computed in \cite{BE3}
as
$[\gamma] = [id] \oplus [\a^+_{1,0}\a^-_{1,1}]$.

\noindent (ii) {\it Example:} ${\cE^{(12)}}$, $\SUd_9\subset \mathit(E_6)_1$.

This modular invariant is at level 12:

\[
\cZ_{\cE^{(12)}} = |\chi_{0,0} + \chi_{9,0} +
\chi_{0,9} + \chi_{4,1} + \chi_{1,4} + \chi_{ 4,4}|^2
+ 2 |\chi_{2,2} + \chi_{5,2} + \chi_{2,5}|^2
\,.
\]

\noindent It is obtained from the conformal embedding
$\SUd_9\subset \mathit(E_6)_1$,
with dual canonical endomorphism given by the vacuum sector:

\be
\label{cansud}
[\theta] =
[ \la_{0,0}] + [\la_{9,0}] +
[\la_{0,9}] + [\la_{4,1}] +[ \la_{1,4}] +[ \la_{ 4,4}].
\ee
This modular invariant can also be realized from the dual canonical
endomorphism
\[
\oplus_\la Z_{\la,\co\la} [\la] = [\la_{0,0}] + [\la_{9,0}] +
[\la_{0,9}] + [\la_{4,1}] + [\la_{1,4}] + [\la_{ 4,4}]
+ 2[\la_{2,2}] + 2[\la_{5,2}] + 2[\la_{2,5}] ,
\]
\noindent where the sum is over all sectors in $\NXN$
using \cite {BER2} the extension $N \subset M \rtimes\bbZ_3$, as $\mathit E_6$
at level $1$ has $\bbZ_3$ fusion rules.

Now the canonical endomorphism corresponding to
\erf{cansud} was computed in \cite{BE3}
as
\[
[\gamma] = [id] \oplus  [\a^+_{1,0}\a^-_{1,1}]
\oplus  [\a^+_{1,1}\a^-_{1,0}]
\oplus  [\a^+_{2,0}\a^-_{2,2}]
\oplus  [\a^+_{2,2}\a^-_{2,0}]
\oplus  [\a^+_{2,1}\a^-_{2,1}].
\]

So we expect six sheets in the full $M$-$M$ system, but this
is where a surprise appears. We do not get six copies
of the $N$-$M$ graph $\cE^{(12)}_1$.  We only get three copies
of $\cE^{(12)}_1$, located at the three sectors
$[\id]$, $[\a^+_{1,0}\a^-_{1,1}]$,
$ [\a^+_{1,1}\a^-_{1,0}]$
in the $\MXM$ graph and three copies of the
isospectral graph $\cE^{(12)}_2$  located at the three
sectors
$[\a^+_{2,0}\a^-_{2,2}]$,
$ [\a^+_{2,2}\a^-_{2,0}]$,
$ [\a^+_{2,1}\a^-_{2,1}]$
 in $\MXM$.

Here we  show that for the conformal inclusion
$\SUd_9\subset(\rmE_6)_1$, for which we have
six $\MXMp$ orbits in $\MXM$, we find three copies of
the graph $\cE^{(12)}_1$ and three copies of $\cE^{(12)}_2$.

Let us  draw the fusion graph of the generator $[\a^+_{(1,0)}]$
in $\MXM$ in blue.
(We  use the labelling as in \cite[Fig.\ 12]{BE3}.)
The vacuum column forces its identity component,
i.e.\ the chiral fusion graph of $[\a^+_{(1,0)}]$,
to be $\cE^{(12)}_1$, see \cite{BE3}.
Now let us think of the fusion graph of $[\a^-_{(1,0)}]$
in $\MXM$ as being red. We now use the fact
that $\cE^{(12)}_j$, $j=1,2,3$, exhaust the list
of isospectral graphs. The connected components of the
red graph will correspond to nimreps and hence
must be $\cE^{(12)}_j$, $j=1,2,3$. (Note that
the modular invariant obeys $Z^*Z=6Z$, hence we
must have six layers.) Which one of the three graphs can
touch the vertices of the blue $\cE^{(12)}_1$?
At the identity vertex this is clearly the other chiral graph,
determined by the vacuum row to be (a red) $\cE^{(12)}_1$.
These two (blue and red) $\cE^{(12)}_1$'s intersect exactly
on the marked (ambichiral) vertices. The other red ``coset''
graphs will connect the other $\MXMo$ fusion orbits in
$\MXMp$. Now the $\MXMo$ fusion orbits are just the
$\bbZ_3$ symmetry orbits of $\cE^{(12)}_1$.
Thus we will have six red layers: The first is the already
determined $\cE^{(12)}_1$ corresponding to the $\MXMo$
orbit of $\id$. Then there will be one layer
connecting $[\a^+_{(1,0)}]$, $[\a^{+(1)}_{(3,1)}]$ and
$[\a^{+(2)}_{(3,1)}]$, similarly one layer connecting
$[\a^+_{(1,1)}]$, $[\a^{+(1)}_{(3,2)}]$ and
$[\a^{+(2)}_{(3,2)}]$, and finally each $\MXMo$ fixed
point $[\a^+_{(2,0)}]$, $[\a^+_{(2,1)}]$ and
$[\a^+_{(2,2)}]$ are connected to one red layer.
To determine the red layer which touches $[\a^+_{(1,0)}]$,
we compute
\[
\langle\a^+_{(1,0)}\a^-_{(1,0)},\a^+_{(1,0)}\a^-_{(1,0)}\rangle=
\langle\a^+_{(1,0)}\a^+_{(1,1)},\a^-_{(1,0)}\a^-_{(1,1)}\rangle=1.
\]
Thus $[\a^+_{(1,0)}]$ has only one target vertex on the red graph.
Hence we must have here either one of the three extremal
vertices of $\cE^{(12)}_1$ or the unique isolated extremal
vertex of $\cE^{(12)}_2$. Since $\cE^{(12)}_3$ does not have
such a vertex, this one is ruled out here.
Now note that the target vertices of these extremal vertices
have itself two and four target vertices for $\cE^{(12)}_1$
and $\cE^{(12)}_2$, respectively.
But since
\[
\langle\a^+_{(1,0)}\a^-_{(1,0)}\a^-_{(1,0)},
\a^+_{(1,0)}\a^-_{(1,0)}\a^-_{(1,0)}\rangle=2
\]
we conclude that a red $\cE^{(12)}_1$ touches $[\a^+_{(1,0)}]$.
The same is checked for $[\a^+_{(1,1)}]$, and it cannot lie on
the same red $\cE^{(12)}_1$ as $[\a^+_{(1,0)}]$ since this
would mean that one is the fusion product of the other by an
ambichiral sector.
Next we check what red graph touches $[\a^+_{(2,0)}]$.
Since
\[
\langle\a^+_{(2,0)}\a^-_{(1,0)},\a^+_{(2,0)}\a^-_{(1,0)}\rangle=1
\]
we must again locate an extremal vertex of $\cE^{(12)}_1$
or $\cE^{(12)}_2$ here. But now
\[
\langle\a^+_{(2,0)}\a^-_{(1,0)}\a^-_{(1,0)},
\a^+_{(2,0)}\a^-_{(1,0)}\a^-_{(1,0)}\rangle=4
\]
forces us to select $\cE^{(12)}_2$.
(We used $[\a^+_{(2,0)}][\a^+_{(2,2)}]
=[\id]\oplus[\a^+_{(2,1)}]\oplus[\a^+_{(4,2)}]$.)
A similar argument applies
to $[\a^+_{(2,1)}]$ and $[\a^+_{(2,2)}]$.
Thus we have indeed found three layers of $\cE^{(12)}_1$
and three layers of $\cE^{(12)}_2$.


Di Francesco and Zuber actually produced
three isospectral graphs $\cE^{(12)}_i$, $i$ = 1,2,3,
whose spectrum
reproduced the diagonal part of the modular invariant
${\cE^{(12)}}$ ($\SUd_9\subset \mathit(E_6)_1$),
 and we realized two of those graphs $\cE^{(12)}_1$
and $\cE^{(12)}_2$ in
\cite{BER1} . The third was apparently eliminated by
\cite{O1}.
We certainly know that $\cE^{(12)}_3$ does not appear
in a ``natural'' way in the sense that we have some
subfactor $N\subset M$ producing $\cE^{(12)}_3$
as $M$-$N$ graph in the following sense:
We know that such a subfactor would have intermediate
subfactors $N\subset M_+=M_-$ producing the same invariant
$Z_{\cE^{(12)}}$ and with $M_+$-$N$ graph $\cE^{(12)}_1$.
This subfactor could not have the ``natural'' property
that the dual canonical endomorphism of  $M_+\subset M$
decomposes exclusively into ambichiral sectors.
This is because we know that the only irreducible
braided extensions (relative to the ambichiral system)
are the trivial one $M_+\subset M=M_+$ and
$M_+\subset M=M_+\rtimes\bbZ_3$ where in turn
$N\subset M$ produces $\cE^{(12)}_1$ and $\cE^{(12)}_2$,
respectively \cite{BER1,BER2}.

\noindent (iii) {\it Example:} ${\cE^{(24)}}$
$\SUd_{21}\subset(\mathrm{E}_7)_1$:

The corresponding modular invariant reads
\[ \begin{array}{ll}
Z_{\mathcal{E}^{(24)}}
=& | \chi_{0,0} + \chi_{{21},0} + \chi_{{21},{21}} +
\chi_{8,4} + \chi_{{17},4} + \chi_{{17},{13}} \\[.4em]
& \qquad + \chi_{{11},1} + \chi_{{11},{10}} + \chi_{{20},{10}} +
\chi_{{12},6} + \chi_{{15},6} + \chi_{{15},9} |^2 \\[.4em]
& + | \chi_{6,0} + \chi_{{21},6} + \chi_{{15},{15}} +
\chi_{{15},0} + \chi_{{21},{15}} + \chi_{6,6} \\[.4em]
& \qquad + \chi_{{11},4} + \chi_{{17},7} + \chi_{{14},{10}} +
\chi_{{11},7} + \chi_{{14},4} + \chi_{{17},{10}} |^2 \,,
\end{array} \]
therefore
\[ \begin{array}{ll}
[\theta] =& [\la_{0,0}] \oplus [\la_{{21},0}] \oplus [\la_{{21},{21}}]
 \oplus [\la_{8,4}]  \oplus [\la_{{17},4}]  \oplus [\la_{{17},{13}}] \\[.4em]
& \qquad [\la_{{11},1}] \oplus [\la_{{11},{10}}] \oplus [\la_{{20},{10}}]
 \oplus [\la_{{12},6}] \oplus [\la_{{15},6}] \oplus [\la_{{15},9}]  \,.
\end{array} \]

\noindent Taking the extension $N \subset M \rtimes\bbZ_2$ \cite{BER2},
 as the extended
system $\mathit E_7$
at level $1$ has $\bbZ_2$ fusion rules, the modular invariant can also be
realised from $\oplus_\la Z_{\la,\co\la} [\la] $
 where the sum is over all sectors in $\NXN$.

\subsection{Towards a pattern}
\labl{pattern}

We have seen that we have exactly
$\sum_\la Z_{0,\la}^2$ (respectively $\sum_\la Z_{\la,0}^2$)
$\MXMp$ (respectively $\MXMm$) orbits in $\MXM$.
These intersect with $\MXMm$ (respectively $\MXMp$),
i.e.\ with the $\MXMm$ (respectively $\MXMp$) orbit
containing $[\id]$, precisely on its $\MXMo$ orbits.
We are interested in the particular shape of the $\MXMp$
or $\MXMm$ orbits in the full system $\MXM$.
For all examples we know, the products $ZZ^*$
and $Z^*Z$ are integral linear combinations of
physical invariants, and the linear combination
corresponds precisely to the decomposition of the
full system in chiral orbits.
Note that each $\MXMpm$ orbit must be a nimrep.
As long as we have a one-to-one correspondence
between irreducible\footnote{We  do not  mean irreducibilty in
the usual sense for representations here --- this would mean
``one-dimensional'' since our braided systems $\NXN$ are
commutative. Here we rather mean irreducibility in the
sense that the sum of the representation matrices is
irreducible (in the sense of \cite{GHJ}).} nimreps
and diagonals of modular invariants we find that
at least the diagonal part of $ZZ^*$ and $Z^*Z$ can be
written as a positive integral linear combination of
diagonal parts of modular invariants.
Since there are no 
distinct\footnote{Here we do not worry about the 
distinction between a sufferable  modular
invariant and its transpose, which can be obtained from the
same subfactor by reversing the braiding}
modular invariants known
sharing the same diagonal part,
this is a strong indication that there is indeed
a general rule.

We can write $Z$ in terms of rectangular branching matrices
as $Z=B_+^* B_-$ so that $ZZ^*=B_+^*B_-B_-^*B_+$
and $Z^*Z=B_-^*B_+B_+^*B_-$.
Let us look at the sandwiched $B_\pm B_\pm^*$ which
must be invariant under the extended
$S$- and
$T$-matrices
thanks to the intertwining rules of \cite[Thm.\ 6.5]{BE4}.
The extended $S$- and $T$-matrices have at most permutation
invariants. If these invariants in fact span the entire
commutant of $S$ and $T$ (may well be in general)
 then $B_\pm B_\pm^*$ must be a linear
combination of these permutations. Unfortunately, it is
not clear whether this is always an \textit{integral}
linear combination.

It is very instructive to look at some examples.
Even \typei\ invariants are interesting here, i.e.\
when we have $B_+=B_-$. For instance for the
D$_{10}$ invariant of $\SUz_{16}$ we have
\[
B_+B_+^*=\left( \begin{array}{cccccc}
2&0&0&0&0&0\\0&2&0&0&0&0\\0&0&2&0&0&0\\
0&0&0&2&0&0\\0&0&0&0&1&1\\0&0&0&0&1&1
\end{array} \right) = 2\cdot\bfe_4 \oplus
\left( \begin{array}{cc}
1&1\\1&1 \end{array} \right)
=\bfe_6 + t_0,
\]
where $t_0$ is is the transposition matrix which
exchanges the two marked vertices $[\a_8^{(j)}]$, $j=1,2$,
on the short legs of D$_{10}$. For the E$_7$ invariant we
have $B_-=\Pi B_+$, where the permutation $\Pi$ is either
$t_j$, $j=1,2$, the transpositions exchanging $[\a_8^{(j)}]$
with the marked vertex $[\a_2]$, or one of the two cyclic
permutations $c_1$, $c_2$. For example, if $\Pi=t_1$,
then $B_-B_-^*=\bfe+t_2$.
Next let us consider the $\cD^{(12)}$ invariant of $\SUd_9$.
Here we find
\[
B_+B_+^*=3\cdot\bfe_6\oplus\left( \begin{array}{ccc}
1&1&1\\1&1&1\\1&1&1\end{array}\right)=\bfe_9+c_1+c_2,
\]
where here $c_1$, $c_2$ denote the two non-trivial
cyclic permutations of the three fixed points
$[\a^{(j)}_{(6,3)}]$, $j=1,2,3$.
For simple current invariants with a single full fixed point
we probably have a sum over all cyclic permutations of
the fixed point constituents in general.

For the conformal inclusion invariant $\cE^{(12)}$
we find
\[
B_+B_+^*=\left( \begin{array}{ccc}
6&0&0\\0&3&3\\0&3&3\end{array}\right)=3\cdot\bfe_3+3\cdot C,
\]
with the $\bbZ_3$ charge conjugation $C$ exchanging the two
non-trivial marked vertices $[\eta_1]$ and $[\eta_2]$.
These numbers reflect exactly the appearance of three times
$\cE^{(12)}_1$ which corresponds to $\bfe$ and three times
$\cE^{(12)}_2$ which corresponds to $C$.
The very special property of this example is that
the orbifold corresponding to charge conjugation changes
the graph non-trivially, $\cE^{(12)}_2$ is the $\bbZ_3$ orbifold
of $\cE^{(12)}_1$ whereas the modular invariant is self-conjugate,
$Z=CZ$. We do not know any other example where this happens.
Other examples for self-conjugate modular invariants which are
non-self-conjugate on the extended level are D$_{4\varrho}$
for $\SUz$. But for D$_4$, the conjugation of the extended
conjugation is obtained by a $\bbZ_3$ orbifold and D$_4$
is its own $\bbZ_3$ orbifold. For $\varrho>1$, the conjugation
will no longer be obtained as an orbifold since we do not
have a simple current group as extended theory, but apparently
the D$_{4\varrho}$ graphs are in general identical with
there non-group-like orbifolds.
Another example is the conformal inclusion
$\SUf_2\subset{\mathit{SU}}(6)_1$, for which the
extended conjugation is also obtained by a $\bbZ_3$ orbifold,
however, the graphs are their own $\bbZ_3$ orbifolds.

Another strange but different case is the conformal inclusion
$\SUf_6\subset{\mathit{SU}}(10)_1$ invariant for which
\[
B_+B_+^*=\left( \begin{array}{cc}
4&0\\0&4\end{array}\right)\oplus
\left( \begin{array}{cc}
3&1\\1&3\end{array}\right)\otimes\bfe_4=3\cdot\bfe_{10}+C,
\]
with the $\bbZ_{10}$ charge conjugation $C$. An in fact,
here we expect three layers of the chiral graph and one
layer of the conjugation graph in the entire $M$-$M$ system.
(See subsection below.)

\subsection{Interesting invariants of $\SUn_n$}
\labl{interesting}

In the conclusions of \cite{BE4} we speculated about
modular invariants which look like \typei\ or \typeii\
but really come from heterotic extensions, i.e.\
for which we have different intermediate local subfactors
$M_+\neq M_-$. By the results of \cite[Sect.\ 4]{BE4}
this means that at least for one $\la$ we have
$\Hom(\id,\a^+_\la)\neq\Hom(\id,\a^-_\la)$
in spite of $Z_{\la,0}=Z_{0,\la}$. Since
$\Hom(\id,\a^\pm_\la)\subset\Hom(\iota,\iota\la)$
this will necessarily require $\langle\canr,\la\rangle\ge2$
for such $\la$. In \cite{BE4}, we pointed out that such
a case may be possible but also that did not  know
of an example.

Here are  examples, actually making use
of the heterotic ${\mathit{SO}}(16\ell)_1$ modular
invariants ($\ell=1,2,3,...$) treated in \cite{BE4}.
For this we consider once more the series of
conformal inclusions $\SUn_n\subset{\mathit{SO}}(n^2-1)_1$.
Note that for $n=7,9,15,17,23,...$, i.e.\ for
$n=8r\pm1$, $r=1,2,3,...$, the number $n^2-1$ is
a multiple of 16, so that the ambient algebra
has a heterotic extension. (The simplest case is therefore
${\mathit{SU}}(7)_7\subset{\mathit{SO}}(48)_1$.)
Using the standard labelling for the sectors of
${\mathit{SO}}(16\ell)_1$, the two heterotic
invariants can be written as
\[
\cZ=\chi^0(\chi^0)^* + \chi^\rms (\chi^0)^*
 +\chi^0(\chi^\rmc)^* + \chi^\rms (\chi^\rmc)^*
\]
and $\cZ^*$ (their coupling matrices are denoted
by $Q$ and $\tmat Q$ in \cite{BE4}).
Now let $N\subset\tilde M$ denote the conformal inclusion
subfactor of $\SUn_n\subset{\mathit{SO}}(n^2-1)_1$
for $n=8r\pm1$, $r=1,2,3,...$.
As we know from \cite[Sect.\ 7]{BE4}, there is a crossed
product extension by all ${\mathit{SO}}(64r^2\pm16r)_1$ sectors
v, s, c (and $0$)
$\tilde M \subset M= \tilde M \rtimes (\bbZ_2\times\bbZ_2)$
producing $\cZ$ and $\cZ^*$ (using braiding and its opposite).
The local intermediate extensions are different, namely
the $\bbZ_2$ extensions corresponding to s and c separately.
So let us consider the subfactor $N\subset M$.
Then its maximal local intermediate extensions will
therefore be
$M_+=\tilde M \rtimes_\rms \bbZ_2$ and
$M_-=\tilde M \rtimes_\rmc \bbZ_2$.
Nevertheless the ${\mathit{SU}}(8r\pm1)_{8r\pm1}$ invariant
arising from $N\subset M$ does not seem to have non-symmetric
vacuum coupling ---  all known  $\SUn$ invariants
are symmetric.
Therefore we expect that the sectors s and c of
${\mathit{SO}}(64r^2\pm16r)_1$ will have the
same branching rules, i.e.\ have the same restriction
to ${\mathit{SU}}(8r\pm1)_{8r\pm1}$.
(This is quite natural due to the similarity of the
sectors s and c of ${\mathit{SO}}(n)_1$.
For $n$ an odd multiple of $8$ this similarity even
covers the sector v and e.g.\ for the conformal
inclusion $\SUd_3\subset{\mathit{SO}}(8)_1$ all three
sectors v, s, c have the same $\SUd$ restriction.)
Anyway, then here we have a heterotic extension, but the
identical branching rules of s and c would imply that $\cZ$,
when written in ${\mathit{SU}}(8r\pm1)_{8r\pm1}$ characters,
has symmetric coupling matrix and looks and in particular
does not look heterotic anymore. In fact, upon restriction
to ${\mathit{SU}}(8r\pm1)_{8r\pm1}$, both $\cZ$ and $\cZ^*$
will be identical with the invariants $|\chi^0+\chi^\rms|^2$
and $|\chi^0+\chi^\rmc|^2$.
Hence there will be a 4-fold degeneracy.
Due to the permutation $\rms\leftrightarrow\rmc$,
the original conformal inclusion invariant
$|\chi^0|^2+|\chi^\rmv|^2+|\chi^\rms|^2+|\chi^\rmc|^2$
will be two-fold degenerate.

Now let $\la$ be an ${\mathit{SU}}(8r\pm1)_{8r\pm1}$
sector appearing in the restriction of s and hence of c.
Since $\cZ$ contains $\chi^\rms (\chi^0)^*$ and
$\chi^0(\chi^\rmc)^*$ it follows that
$Z_{\la,0}=Z_{0,\la}$ is non-zero. On the other hand,
since the dual canonical endomorphism sector $[\canr]$ of
the full subfactor $N\subset M$ is the $\sigma$-restriction
of $[\id]\oplus[\rmv]\oplus[\rms]\oplus[\rmc]$ it follows
that $\langle\canr,\la\rangle\ge2$ --- as it must be.

Now let us  concentrate on the simplest example
of this series, the conformal embedding
${\mathit{SU}}(7)_7\subset{\mathit{SO}}(48)_1$,
which already seems to produce quite interesting
${\mathit{SU}}(7)_7$ modular invariants.
The center of the Weyl alcove, the simple current
fixed point $(1,1,1,1,1,1)$ (or $[6,5,4,3,2,1]$ as a
Young frame) appears in the restriction of s and c,
with a multiplicity 4. Indeed  the
 branching rules are \cite[Eq. (5.30)]{FSS}:
\[
\begin{array}{rl}
0 & \longrightarrow \Big( (0,0,0,0,0,0) \oplus
  (1,0,0,2,1,0) \oplus (0,1,0,0,0,2) \\
  & \qquad \oplus (0,2,0,0,2,0) \oplus
  (2,1,0,1,0,1) \Big) \times\bbZ_7 \oplus (1,1,1,1,1,1) \\[.8em]
\rmv & \longrightarrow \Big( (1,0,0,0,0,1) \oplus
  (3,0,0,1,0,0) \oplus (1,2,0,1,1,0) \oplus (1,1,0,0,1,1)\Big)
  \times\bbZ_7 \\[.8em]
\rms,\rmc & \longrightarrow 4 \cdot (1,1,1,1,1,1),
\end{array}
\]
where ``$\times\bbZ_7$'' means that the entire $\bbZ_7$ simple current
orbit has to be taken.  Note that then indeed the modular invariants $\bfe$
or $W$ of ${\mathit{SO}}(48)_1$ (in the notation of \cite[Sect.\ 7]{BE4})
restrict to the same ${\mathit{SU}}(7)_7$ invariant, let us call it
$Z_\bfe$, and similarly a different ${\mathit{SU}}(7)_7$ invariant, let us
call it $Z_\rms$, is obtained from either $X_\rms$, $X_\rmc$, $Q$ or $\tmat
Q$ of ${\mathit{SO}}(48)_1$ (i.e.\ the latter is the specialization of the
above $\cZ$ or $\cZ^*$.)  Also note that $Z_\bfe$, as it arises from the
diagonal invariant $|\chi^0|^2+|\chi^\rmv|^2+|\chi^\rms|^2+|\chi^\rmc|^2$,
has only two (large) non-zero matrix blocks, because the identical s and c
blocks intersect with the vacuum block --- this is actually the first
modular invariant with this property we have encountered so far.  The first
block, including the vacuum, is a $36\times36$ block, containing 1's
everywhere except a single $33=1+2\cdot4^2$ on the corner corresponding to
the label $(1,1,1,1,1,1)$.  Then there is a $28\times28$ block of 1's
coming from v.  The other invariant $Z_\rms$, arising from
$|\chi^0+\chi^\rms|^2$ (either from $X_\rms$, $X_\rmc$, $Q$ or $\tmat Q$ in
the notation of \cite[Sect.\ 7]{BE4}) has only a single $36\times36$ block,
containing a $35\times35$ block of 1's being cornered by a row and a column
of 35 entries 5, and they meet with a 25 on the corner corresponding to the
label $(1,1,1,1,1,1)$.  Anyway, these seem to be interesting modular
invariants.  First note that we have the curious multiplication rules
$Z_\bfe\cdot Z_\bfe=28 Z_\bfe + 8 Z_\rms$ and $Z_\rms\cdot Z_\rms=60
Z_\rms$.  (Clearly both $Z$'s are selfconjugate in both senses, i.e.\
$Z=CZ$ and $Z=Z^*$.)  Since $\tr Z_\bfe=96$ and $\tr Z_\rms=60$ we will
have $\#\MXN=96=\#\MXMpm$, $\#\MXM=3168$ for (the two-fold degenerate)
$Z_\bfe$, and $\#\MXN=60=\#\MXMpm$, $\#\MXM=3600$ for (the 4-fold
degenerate) $Z_\rms$, and that for $Z_\bfe$ the full system will decompose
into 28 copies of its chiral graph plus 8 copies of the chiral graph for
$Z_\rms$ whereas we expect for $Z_\rms$ itself that the full system will
decompose into 60 layers of its own chiral graph.  (Note that
${\mathit{SU}}(7)_7$ has 1716 primaries.)

It is tempting to conjecture that
for \typei\ invariants, this fusion graph will always
consist exclusively of  copies of the chiral graph.
This is however not the case, as for instance for the
$\cE^{(12)}$ modular invariant of $\SUd_9$ the full system
contains besides 3 copies of $\cE^{(12)}_1$ also 3 copies
of the isospectral graph $\cE^{(12)}_2$, see below.
Moreover, even for \typei\ invariance the product
$ZZ^*$ is not necessarily a multiple of $Z$.
For instance the modular invariant arising from
the conformal inclusion $\SUf_6\subset{\mathit{SU}}(10)_1$
fulfills $ZZ^*=3Z+CZ$, see \cite{BER2}.

\vspace{0.2cm}\addtolength{\baselineskip}{-2pt}
\begin{footnotesize}
\noindent{\it Acknowledgement.}
This work was completed during visits to the MSRI programme
on Operator Algebras in 2000-2001.
 I am grateful for the organizers
and MSRI for their invitation and financial support. I would like to thank
Jens B\"ockenhauer for discussions on this work
as well as for our collaboration in 1996--2000, and Karl-Henning 
Rehren for comments on a preliminary version of this manuscript.
\end{footnotesize}
\vspace{0.5cm}



\newcommand\bitem[2]{\bibitem{#1}{#2}}

\def\aam              {Acta Appl.\ Math. }
\def\aip              {Ann.\ Inst.\ H.\ Poincar\'e (Phys.\ Th\'eor.) }
\def\cmp              {Com\-mun.\ Math.\ Phys. }
\def\duke             {Duke Math.\ J. }
\def\ijm              {Intern.\ J. Math. }
\def\jfa              {J.\ Funct.\ Anal. }
\def\jmp              {J.\ Math.\ Phys. }
\def\lmp              {Lett.\ Math.\ Phys. }
\def\rmp              {Rev.\ Math.\ Phys. }
\def\inv              {Invent.\ Math. }
\def\mpl              {Mod.\ Phys.\ Lett. }
\def\nup              {Nucl.\ Phys. }
\def\nupp             {Nucl.\ Phys.\ (Proc.\ Suppl.) }
\def\adma             {Adv.\ Math. }
\def\physa            {Physica {\bf A} }
\def\ijmp             {Int.\ J.\ Mod.\ Phys. }
\def\jp               {J.\ Phys. }
\def\fdp              {Fortschr.\ Phys. }
\def\plb              {Phys.\ Lett.\ {\bf B}}
\def\rims             {Publ.\ RIMS, Kyoto Univ. }


\begin{footnotesize}

\end{footnotesize}
\end{document}